\documentclass[12pt,reqno,intlim]{amsart}

\usepackage{amssymb,amsmath}

\usepackage[arrow,matrix,curve]{xy}

\newdir{ >}{{}*!/-5pt/\dir{>}}

\newcommand{\nc}{\newcommand}

\nc{\bC}{\bold{C}} \nc{\bN}{\Bbb{N}} \nc{\cF}{\mathcal{F}}
\nc{\cE}{\mathcal{E}} \nc{\cR}{\mathcal{R}} \nc{\cM}{\mathcal{M}}
\nc{\al}{\alpha} \nc{\bt}{\beta} \nc{\gm}{\gamma} \nc{\dl}{\delta}
\nc{\om}{\omega} \nc{\sg}{\sigma} \nc{\Sg}{\Sigma} \nc{\vf}{\varphi}
\nc{\ve}{\varepsilon} \nc{\os}{\overset} \nc{\ol}{\overline}
\nc{\ul}{\underline} \nc{\us}{\underset} \nc{\sbs}{\subset}
\nc{\bsl}{\backslash} \nc{\Ra}{\Rightarrow}
\nc{\lra}{\longrightarrow} \nc{\all}{\allowdisplaybreaks}

\nc{\Codes}{\operatorname{{\bold{Codes}}}}
\nc{\RegMono}{\operatorname{\mathcal{R}{\rm{eg}\mathcal{M}{\rm{ono}\!}}}}
\nc{\RegEpi}{\operatorname{\mathcal{R}{\rm{eg}\mathcal{E}{\rm{pi}\!}}}}
\nc{\Mn}{\operatorname{\mathcal{M}{\rm{ono}\!}}}
\nc{\Ep}{\operatorname{\mathcal{E}{\rm{pi}\!}}}
\nc{\Rg}{\operatorname{\mathcal{R}{\rm{eg}\!}}}
\nc{\Ob}{\operatorname{Ob\!}}

\numberwithin{equation}{section}

\newtheorem{theo}{\ \ \ Theorem}[section]
\newtheorem{lem}[theo]{\ \ \ Lemma}

\newtheorem{cor}[theo]{\ \ \ Corollary}

\theoremstyle{definition}

\theoremstyle{remark}

\begin{document}

\title[]
{Effective codescent morphisms of Hausdorff topological spaces}

\author{Dali Zangurashvili}

\maketitle

\begin{abstract}
Our earlier results on effective codescent morphisms of Hausdorff topological spaces are strengthened and complemented. In particular, it is proved that any embedding $p:B\rightarrowtail E$ with compact $B$ and normal $E$ is an effective codescent morphism in the category of Hausdorff topological spaces.  
\end{abstract}
\bigskip



\section{Introduction}

Let $\mathcal{C}$ be a category with pushouts, and $p: B\rightarrow E$ be a morphism in $\mathcal{C}$. Recall that $p$ induces the following adjunction between the coslice categories: the left adjoint is the change-of-cobase functor $$p^{*}: E/\mathcal{C}\rightarrow B/\mathcal{C}$$ (pushing out along $p$), while the right adjoint is the functor $$p^{!}: B/\mathcal{C}\rightarrow E/\mathcal{C}$$ that composes with $p$ on the right. The morphism $p$ is called a codescent morphism (resp. an effective codescent morphism) if the functor $p^{*}$ is precomonadic (resp. comonadic) \cite{JT}. If a category has equalizers, then a morphism is a codescent morphism if and only if it is a universal regular monomorphism, i.e., a regular monomorphism such that its pushout along any morphism is again a regular monomorphism \cite{JT1}. 

In \cite{Z2}, we proved that every codescent morphism is effective in the category $Haus$ of Hausdorff topological spaces. However, not every regular monomorphism (= closed embedding) is a codescent morphism in this category, as it follows from the example provided in \cite[Remark 5.15]{K}. Therefore, the problem of characterizing codescent morphisms in the category $Haus$ arose. In \cite{Z2}, we gave a necessary condition for a regular monomorphism  $B\rightarrowtail E$ of Hausdorff topological spaces to be universal. Moreover, it was proved that if $B$ is compact and $E$ is regular, then this condition is also sufficient.

In the present paper, we strengthen and complement the original results. In particular, the requirement that $E$ is regular is removed from the second assertion. Moreover, we give a new sufficient condition for a closed embedding $B\rightarrowtail E$ of Hausdorff topological spaces to be a universal regular monomorphism. Besides, we prove that if $B$ is compact, then this condition is also necessary.  With the aid of these facts, it is shown that if $B$ is compact and $E$ is normal, then $p$ is an effective codescent morphism in the category of Hausdorff topological spaces. In particular, this is the case if $E$ is compact and Hausdorff.

Financial support from  Shota Rustaveli  National Science Foundation of Georgia
(Ref.: FR-24-8249) is gratefully acknowledged.

\section{Effective codescent morphisms}
\vskip+2mm

We begin with the following statement.

\begin{theo} \cite[Example 3.11]{Z2}
Any codescent morphism of Hausdorff topological spaces is effective.
\end{theo}

Let $X$ be a topological space. Recall that subsets $Y$ and $Z$ of $X$ are called completely separable if there exists a continuous mapping $f$ from $X$ to the real number interval $[0,1]$ such that $f(Y)=0$  and $f(Z)=1$.

\begin{lem}
The closures of completely separable subsets $Y$ and $Z$ of a topological space $X$ are disjoint. The converse is true if $X$ is normal.
\end{lem}

\begin{proof}
One obviously has $Y\subseteq f^{-1}\lbrace 0\rbrace$ and $Z\subseteq f^{-1}\lbrace 1\rbrace$. Since $f^{-1}\lbrace 0\rbrace$ and $f^{-1}\lbrace 1\rbrace$ are closed,
$\overline{Y}\subseteq f^{-1}\lbrace 0\rbrace$ and $\overline{Z}\subseteq f^{-1}\lbrace 1\rbrace$. This implies the first claim. The second claim immediately follows from Urysohn's lemma.
\end{proof}

Before continue, recall the construction of pushouts
\begin{equation}
\xymatrix{
B\ar[r]^{p}\ar[d]_{f}&E\ar[d]^{f'}\\
A\ar[r]^{p'}&D}
\end{equation}
in the category $Top$ of topological spaces. First of all, recall that the forgetful functor from $Top$ to the category of sets is topological. This implies that square (2.1) is a pushout in the category of sets too. Moreover, $D$ is equipped with the final topology:
a subset $W$ of $D$ is open if and only if the inverse images $p'^{-1}(W)$ and  $f'^{-1}(W)$ are open in $A$ and $E$ respectively. 

In the sequel, we will assume that $p$ is injective. Then $D$ can be viewed as the union $A\cup (E\setminus B)$ (without loss of generality, we assume that $A$ and $E\setminus B$ are disjoint) with $p'$ being the embedding, while $f'$ being the mapping defined as follows: 

\begin{center}
$f'(b)=f(b)$ if $b\in B$, and $f(x)=x$ if $x\in E\setminus B$.  
\end{center}

 Obviously, any subset $W$ of $D$ can be uniquely written  as the union $X\cup Y$ where $X\subseteq A$ and $Y\subseteq E\setminus B$. Then we have:
\vskip+1mm 
\noindent \textit{a subset $W$ of $D$ is open if and only if $X$ is open in $A$ and the union $f^{-1}(X)\cup Y$ is open in $E$}.

\begin{lem}
Let (2.1)
be a pushout in the category $Top$, and $p$ be a closed embedding. Let $a_1$ and $a_2$ be distinct points of the space $A$. Assume that there are disjoint open subsets $U_1$ and $U_2$ of $A$ containing resp. $a_1$ and $a_2$,  and also there are disjoint open subsets $V_1$ and $V_2$ of $E$ such that $V_i\cap B=f^{-1}(U_i)$ ($i=1,2$). Then $p'(a_1)$ and $p'(a_2)$ have disjoint neighbourhoods. 
\end{lem}

\begin{proof}
The subsets $W_i=U_i\cup (V_i\setminus B)$ ($i=1,2$) of $D$ are obviously open, disjoint, and contain $p'(a_i)$.
\end{proof}

\begin{lem}
 Let (2.1) be a pushout in the category $Top$ of topological spaces with Hausdorff $A$, $B$, $E$, and let $p$ be a closed embedding. Assume that $x_1$ and $x_2$ are distinct points of $D$. If $B$ is compact or $E$ is regular, then $x_1$ and $x_2$ have disjoint neighbourhoods (in $D$) in any of the following cases:
(a)  $x_1, x_2\in E\setminus B$; (b) $x_1\in A$ and $x_2\in E\setminus B$.
\end{lem}

\begin{proof}
 (a) Let $V_1$ and $V_2$ be disjoint neighbourhoods of $x_1$ and $x_2$ in $E$. The intersections $V_i\cap (E\setminus B)$ obviously contain $x_i$ ($i=1,2$) and are disjoint and open in $D$. 
 
 (b) There are open disjoint subsets $U$ and $V$ of $E$ with $B\subseteq U$ and $x_2\in V$. In the case where $E$ is regular, this is obvious. If $B$ is compact, this follows from, e. g., \cite[Lemma 26.4]{Mu}. Let $W=A\cup (U\setminus B)$. Obviously $f^{-1}(A)\cup (U\setminus B)=U$. Hence $W$ is an open neighbourhood of $x_1$ in $D$. Now it suffices to observe that $V$ is an open neighbourhood of $x_2$ in $D$ and is disjoint with $W$.
\end{proof}

The following lemma is well-known. For the reader's convenience, we give it with a proof.

\begin{lem}
Let (2.1) be a pushout in the category of topological spaces, and let $p$ and $f$ be closed embeddings. If $A$, $B$ and $E$ are Hausdorff, then so is $D$, and the mappings $p'$ and $f'$ are closed embeddings.
\end{lem}

\begin{proof}
For the particular case considered in  this lemma, the above-mentioned construction of a pushout in the category of topological spaces can be simplified: $D$ is the amalgam, i.e., the union of copies of $A$ and $E$ with the intersection being precisely $B$ (for the sake of brevity, we will not distinguish these copies from $A$ and $E$ resp.). Moreover, a subset $W$ of $D$ is open if and only if $W\cap A$ and $W\cap E$ are open in resp. $A$ and $E$.

Let $x_1$ and $x_2$ be distinct points of $D$. There is only one non-trivial case to consider, it is the one where both $x_1$ and $x_2$ belong to $B$. In that case, there are open neighbourhoods $U_1$, $U_2$ of $x_1$ and $x_2$ in $A$, and $V_1$, $V_2$ of $x_2$ in $E$ such that $U_1$ and $U_2$ are disjoint and $V_1$ and $V_2$ are disjoint. Let $X_i=U_i\cap V_i$  $(i=1,2)$. Then there are open subsets $U'_i$ of $A$ and open subsets $V'_i$ of $E$ with $U'_i\cap B=X_i=V'_i\cap B$  $(i=1,2)$. The subsets $W_i=(U_i\cap U'_i)\cup (V_i\cap V'_i)$ of $D$ are obviously disjoint, open and contain $x_i$  $(i=1,2)$.
\end{proof}

Further, recall that the category of Hausdorff topological spaces is reflective in the category of topological spaces \cite[Proposition 2, p.135]{M} (note that this fact also immediately follows from the epireflective subcategory theorem). This implies the construction of pushouts in the category of Hausdorff topological spaces: it is the outer quadrangle in the following diagram, where the internal square is a pushout in the category of topological spaces, and $\eta$ is the unit of the reflection:

\begin{equation}
\xymatrix{
B\ar[r]^{p}\ar[d]_{f}&E\ar[d]^{f'}\ar[ddr]\\
A\ar[r]^{p'}\ar[drr]_{p''}&D\ar[dr]^{\eta_D \; \; \; \; \; \; \; \; \; \; \; \; \;\; \; \; \; \; \; \; \; \; \; \; \; \;\; \; \; \; \; \; \; \; \; \; \; \; \; \; \;\; \; \; \; \; \; \; \; \; \; \; \; \;\; \; \; \; \; \; \; \; \; \; \; \; \; }\\
&&D'}
\end{equation}

\begin{theo}
For a closed embedding $p:B\rightarrowtail E$ of Hausdorff topological spaces, and the following conditions, one has the implications (iv)$\Rightarrow$(i)$\Rightarrow$(ii)$\Leftarrow$(iii). If $E$ is regular, then also (iii)$\Rightarrow$(iv). If $B$ is compact, then (i), (ii) and (iv) are equivalent:
\vskip+1mm
(i) $p$ is a universal regular monomorphism;
\vskip+2mm
(ii) for any completely separable open subsets $U_1$ and $U_2$ of $B$, there are disjoint open subsets $V_1$ and $V_2$ of $E$ such that $V_i\cap B=U_i$ ($i=1,2$);
\vskip+3mm
(iii) for any disjoint open subsets $U_1$ and $U_2$ of $B$, there are disjoint open subsets $V_1$ and $V_2$ of $E$ such that $V_i\cap B=U_i$ ($i=1,2$);
\vskip+1mm
(iv) for any continuous mapping $f:B\rightarrow A$ with Hausdorff space $A$, the space $D$  in the pushout
\begin{equation}
\xymatrix{
B\ar[r]^{p}\ar[d]_{f}&E\ar[d]^{f'}\\
A\ar[r]^{p'}&D}
\end{equation}
 in the category of topological spaces, is Hausdorff.
\end{theo}

\begin{proof} (i)$\Rightarrow$(ii): Let $U_1$ and $U_2$ be completely separable open subsets of $B$, and let $f:B\rightarrow [0,1]$ be a continuous mapping with $f(U_1)=0$ and $f(U_2)=1$. Consider the  commutative diagram (2.2) with $A$ being the interval $[0.1]$,
where the internal square is a pushout in the category $Top$. As we already mentioned above, the outer quadrangle is a pushout in the category of Hausdorff topological spaces. Therefore, $p''$ is injective. This implies that  $p''(0)$ and $ p''(1)$ have disjoint open neighbourhoods $W'_1$ and $W'_2$. Then  $\eta_D^{-1}(W'_1)$ and $\eta_D^{-1}(W'_2)$ are disjoint open neighbourhoods of $p'(0)$ and $p'(1)$. Let $\eta_D^{-1}(W'_i)=X_i\cup Y_i$ with $X_i\subseteq [0,1]$ and $Y_i\subseteq E\setminus B$. The set  $f^{-1}(X_i)\cup Y_i$ is open and contains $U_i$ ($i=1,2)$.

Since $p$ is an embedding, there are open subsets $O_i$ of $E$ with $U_i=O_i\cap B$ ($i=1,2$). Now it is clear that $V_i=O_i\cap (f^{-1}(X_i)\cup Y_i)$ are the sought-for open subsets of $E$.

The implication (iii)$\Rightarrow$(ii) is obvious.

The implication (iii)$\Rightarrow$(iv) follows from Lemmas 2.3 and 2.4 provided that $E$ is regular.

(iv)$\Rightarrow$(i): It is well-known that the class of embeddings is closed under pushout in the category $Top$ (in other words, $Top$ is a coregular category). Hence $p'$ in diagram (2.3) is an embedding. Applying the construction of pushouts provided at the beginning of this section, it is easy to see that $p'$ is a closed embedding. Now it suffices to observe that (2.3) is a pushout in the category of Hausdorff spaces too.

(ii)$\Rightarrow$(iv): Let $B$ be compact. Consider diagram (2.3) with Hausdorff $A$. First consider the particular case where $A$ is normal. Let $a_1$ and $a_2$ be distinct points of $A$. Since $A$ is normal,  there are disjoint open neighbouhoods $U_1$ and $U_2$ of resp. $a_1$ and $a_2$ such that their closures $\overline{U_1}$ and $\overline{U_2}$ are disjoint as it follows from \cite[Lemma 31.1(b)]{Mu}). By Urysohn's lemma, there is a continuous mapping $f:A\rightarrow [0,1]$ with $f(\overline{U_1})=0$ and $f(\overline{U_2})=1$. This implies that $U_1$ and $U_2$ are completely separable. Then there are disjoint open subsets $V_1$ and $V_2$ of $E$ with $V_i\cap B=U_i$ ($i=1,2$). Therefore, by Lemma 2.3, $a_1$ and $a_2$ have disjoint open neighbourhoods in $D$. Hence Lemma 2.4 implies that the space $D$ in diagram (2.3) is Hausdorff. Also note that the  arguments provided in the proof of (iv)$\Rightarrow$(i) imply that $p'$ is a closed embedding.

Consider now the the case of an arbitrary Hausdorff space $A$. Factorize the mapping $f$ as $f=em$ with a dense continuous mapping $e$ and a closed embedding $m$. Then pushout (2.3) can be represented as a concatenation of the following two pushouts (in $Top$):
\begin{equation}
\xymatrix{
B\ar@{ >->}[r]^{p}\ar[d]_{e}&E\ar[d]\\
C\ar@{ >->}[r]^{g}\ar@{ >->}[d]_{m}&F\ar[d]\\
A\ar[r]^{p'}&D}
\end{equation}

Since $B$ is compact, so is $f(B)$. This implies that $f(B)$ is closed. Therefore, $C=f(B)$, and hence it is normal. According to what already was proved, $F$ is Hausdorff, and $g$ is a closed embedding. Lemma 2.5 implies that $D$ is Hausdorff. 
\end{proof}

 \begin{cor}
 Let $B$ be a clopen (i.e., closed and open) subspace of a topological space $E$. If $B$ is compact or $E$ is regular, then the embedding $B\rightarrowtail E$ is an effective codescent morphism in the category of Hausdorff topological spaces.
 \end{cor}
 
 \begin{proof}
 It suffices to observe that, under the assumptions of this corollary, the condition (iii) is obviously satisfied. 
 \end{proof}
 
 \begin{lem}
 Let $p:B\rightarrowtail E$ be a closed embedding. If $E$ is normal, then the condition (ii) of Theorem  2.6 is satisfied.
 \end{lem}
 
 \begin{proof}
 Let $U_1$ and $U_2$ be open completely separable subsets of $B$.  Then, by Lemma 2.2, their closures $C_1=\overline{U_1}$ and $C_2=\overline{U_2}$ in $B$ are disjoint. Since $B$ is closed in $E$, $C_1$ and $C_2$ are closed in $E$ too. Since $E$ is normal, there are disjoint open subsets $V_1$ and $V_2$ of $E$ containing resp. $C_1$ and $C_2$.
 
 Further, since $p$ is an embedding, there are open subsets $V_1'$ and $V_2'$ of $E$ with $V'_i\cap B=U_i$ ($i=1,2$). Then the intersections $V_i\cap V'_i$ ($i=1,2$) are the sought-for sets. Indeed, they are disjoint, open and  obviously $(V_i\cap V'_i)\cap B=U_i$.
 \end{proof}

Theorems 2.1, 2.6 and Lemma 2.8 immediately imply the following statement. 
 \begin{cor}
Let $p:B\rightarrowtail E$ be a closed embedding. If $B$ is compact and $E$ is normal, then $p$ is an effective codescent morphism in the category of Hausdorff topological spaces.
 \end{cor}

Finally, Corollary 2.9 implies the following assertion.
 
 \begin{cor}
Let $p:B\rightarrowtail E$ be a closed embedding. If $E$ is compact and Hausdorff, then  $p$ is an effective codescent morphism in the category of Hausdorff topological spaces.
\end{cor}

\vskip+2mm
 
\textit{Author's Address:}
\vskip+1mm 

\textit{A. Razmadze Mathematical Institute of}

\textit{Iv. Javakhishvili Tbilisi State University,}

\textit{2 Alexidze Str., Lane II,  0193, Georgia;}
\vskip+1mm
\textit{e-mail: dali.zangurashvili@tsu.ge}

\end{document}